\documentclass[conference]{IEEEtran}
\makeatletter
\def\ps@headings{%
\def\@oddhead{\mbox{}\scriptsize\rightmark \hfil \thepage}%
\def\@evenhead{\scriptsize\thepage \hfil \leftmark\mbox{}}%
\def\@oddfoot{}%
\def\@evenfoot{}}
\makeatother
\usepackage[scriptsize,it]{caption}
\usepackage{amssymb}
\usepackage{amsfonts}
\usepackage{amsmath,amssymb}%
\ifCLASSINFOpdf 
\usepackage[pdftex]{graphicx} 
\usepackage[pdftex,colorlinks,
           bookmarks=true,%
           pdftitle=Using\ Poisson\ processes\ to\ model\ lattice\ cellular\ networks',%
           pdfauthor=B.\ Blaszczyszyn%
           \ -\  M.\ Karray\ -\ P.\ Keeler,]{hyperref}  
\usepackage[numbers,sort&compress]{natbib} 
\usepackage{hypernat} 
\hypersetup{
   bookmarksnumbered,
   pdfstartview={FitH},
   citecolor={blue},
   linkcolor={red},
   urlcolor=[rgb]{0,0.55,0},
   pdfpagemode={UseOutlines}}
\else

 \usepackage[dvips]{graphicx} 
 \usepackage[numbers,sort&compress]{natbib} 
\fi

\newtheorem{thm}{Theorem}
\newenvironment{theorem}{\bf\begin{thm}\rm\em}{\end{thm}} 
\newtheorem{cor}[thm]{Corollary}
\newenvironment{corollary}{\bf\begin{cor}\rm\em}{\end{cor}} 
\newtheorem{lem}[thm]{Lemma}
\newenvironment{lemma}{\bf\begin{lem}\rm\em}{\end{lem}} 
\newtheorem{prop}[thm]{Proposition}
\newenvironment{proposition}{\bf\begin{prop}\rm\em}{\end{prop}} 
\newtheorem{rem}[thm]{Remark}
\newenvironment{remark}{\bf\begin{rem}\rm}{\end{rem}} 

\newcommand{\E}{\mathbf{E}}
\newcommand{\calE}{\mathcal{E}}
\newcommand{\calS}{\mathcal{S}}
\newcommand{\calR}{\mathcal{R}}
\newcommand{\calL}{\mathcal{L}}

\newcommand{\Prob}{\mathbf{P}}
\renewcommand{\P}{\mathbf{P}}
\newcommand{\Ind}{\mathbf{1}}
\newcommand{\R}{\mathbb{R}}

\newcommand{\calN}{\mathcal{N}}
\newcommand{\dB}{\mathrm{dB}}
\newcommand{\logsd}{\sigma_{\dB}}

\begin{document}

\title{Using Poisson processes to model lattice cellular networks}
\author{\IEEEauthorblockN{Bart{\l }omiej~B{\l}aszczyszyn\IEEEauthorrefmark{1},
Mohamed Kadhem
Karray\IEEEauthorrefmark{2} 
and Holger Paul Keeler\IEEEauthorrefmark{1}}
\IEEEauthorblockA{\IEEEauthorrefmark{1}INRIA-ENS,
23 Avenue d'Italie, 75214 
Paris, France
Email: Bartek.Blaszczyszyn@ens.fr, Holger-Paul.Keeler@inria.fr\\
\IEEEauthorrefmark{2}Orange Labs;
38/40 rue G\'{e}n\'{e}ral Leclerc, 92794  
Issy-les-Moulineaux, France
Email: mohamed.karray@orange.com}}
\maketitle

\begin{abstract}
An almost ubiquitous assumption made in the stochastic-analytic approach to study of the quality of user-service in cellular networks is Poisson distribution of base stations, often completed by some specific assumption regarding the distribution of the fading (e.g. Rayleigh).  The former (Poisson) assumption is usually (vaguely) justified in the context of cellular networks, by various irregularities in the real placement of base stations, which ideally should form a lattice (e.g. hexagonal) pattern.  In the first part of this paper we provide a different and rigorous argument justifying the Poisson assumption under sufficiently strong log-normal shadowing observed in the network, in the evaluation of a natural class of the typical-user service-characteristics (including path-loss, interference, signal-to-interference ratio, spectral efficiency). Namely, we present a Poisson-convergence result for a broad range of stationary (including lattice) networks subject to log-normal shadowing of increasing variance.  We show also for the Poisson model that the distribution of all these typical-user service characteristics does not depend on the particular form of the additional fading distribution.  Our approach involves a mapping of 2D network model to 1D image of it ``perceived'' by the typical user. For this image we prove our Poisson convergence result and the invariance of the Poisson limit with respect to the distribution of the additional shadowing or fading.  Moreover, in the second part of the paper we present some new results for Poisson model allowing one to calculate the distribution function of the SINR in its whole domain.  We use them to study and optimize the mean energy efficiency in cellular networks.
\end{abstract}


\section{Introduction}
Cellular networks are being extensively deployed and upgraded in order
to cope with the steady rise of user-traffic. This has created the
need for new and robust analytic techniques to study the quality of
user-service. The ability to tractably model and calculate the quality of user-service
related to 
the signal-to-interference-and-noise ratio (SINR, or SIR when the
noise is neglected)  
will serve as the motivating force behind the work presented here. 

In order to derive analytic techniques, various mathematical
models have been proposed. A common and simplifying model assumption
is that the base stations are located according to a Poisson process
in the plane. In the first section of this paper we recall 
such a model with shadowing and/or fading, and present a simple yet
very useful result, upon which the bulk of the remaining work here
hinges. This result involves mapping the point process on $\R^2$ (that models
the locations of base stations) through the 
distance-loss function and the shadowing and/or fading variables (that model the
propagation losses between these stations and a typical user),
to a {\em point process of propagation losses} on $\R^+$, which are
experienced by this user.  It  allows one to study  
all typical-user characteristics,
 which can be expressed in terms of its  propagation losses
(or received powers, 
e.g.  SINR, etc).
  We observe that 
{\em this new process also forms a Poisson process}, regardless of the
 shadowing and/or fading distribution, and is characterized only by the
moment of this distribution of order $2/\beta$, where $\beta$ is the
distance-loss exponent. 

In the context of an actual deployment of cellular networks, 
lattice (e.g. hexagonal) models for the base station placement
are usually thought of as more pertinent. However, perfect 
lattice models do not seem to allow analytic techniques
for the study of the SINR-based characteristics. 
Hence, the Poisson model is used  and  justified by positioning 
``irregularities'' of the network. It is also considered
as a ``worst-case'' scenario due to its complete randomness
property. Although the validity of these arguments alone may be
questioned, {\em we seek to support the Poisson assumption with a
powerfully new convergence result  when there is sufficiently
strong log-normal shadowing in the network}. 
Namely, provided a network is represented by a sufficiently
large homogeneous  point pattern (whose definition will
later be made precise, and includes general
lattices and various ``perturbed lattice'' models), 
we show that as the variance of
log-normal shadowing increases, the resulting 
propagation losses between the stations and the typical user 
form a stochastic process that converges to the aforementioned 
non-homogeneous Poisson process on $\R$, due to the Poisson model.
In other words, the actual (large but not necessarily Poisson) network 
is perceived by a typical user as an equivalent (infinite) Poisson network, 
provided shadowing  is strong enough, of logarithmic standard deviation
greater than approximately 10dB,
as shown by numerical evidence. This is a  realistic assumption
for outdoor and indoor wireless communications in many urban scenarios.

This result rigorously justifies using a Poisson point
process to model typical-user characteristics
in a wide range of network models, which includes the
hexagonal model and a large class of perturbed lattice models, thus
adding theoretical weight to the work being done under the Poisson
assumption.  

After stating this convergence result, we derive important complementary analytic tools for the
study of the distribution of the SINR in the Poisson model,
which remain valid for any distribution of the
shadowing and/or fading random variables.
We use them to investigate the spectral and energy efficiency in cellular
networks. In particular, we evaluate the  mean  energy efficiency
as a function of the base station transmit power, which allows one to 
optimally tune this latter power.
Poisson and hexagonal networks with and without shadowing are compared
in this context.

\subsubsection*{Related work}
It is beyond the scope of this short introduction to review all 
works with Poisson model of wireless networks.  Some results and
further references may be found in~\cite{FnT,Haenggi2009}. 
More specifically, our results presented in this paper
allow to extend  
to a general fading distribution the explicit expressions for the
distribution function of the SINR derived in~\cite{Dhillon2012}, where 
Rayleigh fading is assumed.
Our representation of the network via the process of propagation
losses is similar tho this  
considered in~\cite{net:Haenggi08tit} for other purposes, mostly to
study the effect of shadowing/fading on connectivity.
Moreover, using the Laplace/Fourier analysis, we manage to characterize   
(no longer explicitly) the SINR over its entire range,
whereas the aforementioned explicit expressions 
are valid only for $\text{SINR}\ge 1$.
The result of Lemma~\ref{l.ModifiedIntensity} appeared in~\cite{shadow2010}.
The infinite Poisson model was statistically 
fitted in~\cite{propag_extabst}
to some real  data
provided by an operator regarding propagation losses
 in order to estimate the parameters of the
propagation loss  
using a simple linear regression model.
The spectral and energy  efficiency in hexagonal networks without shadowing
was studied in~\cite{Karray2010SpectralEfficiency}.
Finally, the convergence result  presented here is
in the spirit of classical limit theorems of point processes, which
are detailed in \cite[Chapter 11]{daleyPPII2008}. In particular, these
theorems show that under specific conditions, the repeated
superposition, thinning or translation of points of a point process
will result in the point process converging in the limit to a Poisson
process.

{\em The remaining part of this paper} is organized as follows. 
In section~\ref{PoissonModel1} we present the basic result on the
network mapping to the process of propagation losses.
The main convergence result is presented in
Section~\ref{s.convergence} and its proof is given Appendix~\ref{App1}.
We revisit the Poisson model in order to study the SINR distribution in Section
\ref{PoissonModel2}.


\section{Infinite Poisson model} \label{PoissonModel1}
For motivation purposes, we first recall the usual cellular network model
based on the Poisson process. 
In particular, we model the geographic locations of the base stations with an homogeneous
Poisson point process $\Phi=\left\{  X_{i}\right\}
_{i\in\mathbb{N}}$ of intensity $\lambda$ on $\mathbb{R}^{2}$, which we refer to as the \emph{infinite
Poisson} model.

We assume that the (typical) user is located at the origin without loss of generality due to the stationarity of $\left\{  X_{i}\right\}  _{i\in\mathbb{N}}$.
Let $l\left(  X_{i}\right)$ be the  \emph{distance loss} between a base station at $X_{i}$\ and the
user, where $l\left(  \cdot\right)$ is given by
\begin{equation}
l\left(  x\right)  =\left(  K\left\vert x\right\vert \right)  ^{\beta}
\label{e.DistanceLoss}%
\end{equation}
for two given positive constants $K$\ and $\beta$. When we incorporate 
\emph{shadowing} (and/or \emph{fading}), the \emph{propagation loss} is
\[
L_{X_{i}}=\frac{l\left(  X_{i}\right)  }{S_{X_{i}}}%
\]
where $\left\{  S_{x}\right\}  _{x\in\mathbb{R}^{2}}$\ is a collection of
independent and identically distributed (iid) positive random
variables. We will sometimes write also $S_{X_i}=S_i$ to simplify the notation. 
Let $S$\ be a random variable having the same distribution as any $S_{x}$.

Note that the {\em power received} at the origin from the station $X_i$,
transmitting with power $P_{X_{i}}$, is
equal to 
\[
p_{X_{i}}=\frac{P_{X_{i}}}{L_{X_{i}}}=\frac{P_{X_{i}}S_{X_{i}}}{l\left(
X_{i}\right)}\,.%
\]
In this paper we do not assume any power control, i.e., $P_{X_i}=P$
for  some given positive constant $P>0$.
In this case (as well in a more general case of iid emitted powers),
including $P_{X_{i}}$ in the associated shadowing random variable, we retrieve
an equivalent model in which the shadowing is $\tilde{S}_{X_{i}}:=P_{X_{i}%
}S_{X_{i}}$\ and the transmitted powers $\tilde{P}_{X_{i}}=1$. 
Henceforth we assume, without loss of generality, that
the transmitted powers are all equal to one, while keeping in mind that the shadowing
random variables now include the effective transmitted powers. This
transformation will slightly simplify our notation. However, when
studying the energy efficiency in Section~\ref{sss.SINR} we will
reintroduce emitted powers to our model.

\subsection{Mapping of the propagation losses of the typical user from $\mathbb{R}^2$ to $\mathbb{R}^+$}
Denote by $\calN=\left\{
L_{X_{i}}\right\}  _{i\in\mathbb{N}}$ the {\em  process of propagation
  losses} experienced by the typical user with respect to the stations
in~$\Phi$. We consider $\calN$ a {\em point process on
$\mathbb{R}^{+}$}. Note that {\em all} characteristics of the typical
user, which can be expressed in terms of its propagation losses
(or received powers, under the aforementioned assumption on
emitted powers, e.g. SIR, SINR, spectral and energy efficiency, etc) 
are determined by the
distribution of $\calN$. This motivates the following simple result
that  appeared, to the best of our knowledge  for the first time,
in~\cite{shadow2010}.
In order to make this presentation more self-contained we present it
with a proof.  
\begin{lemma} 
\label{l.ModifiedIntensity}Assume infinite Poisson model with
distance-loss~(\ref{e.DistanceLoss}) and generic shadowing (and/or
fading) variable satisfying
\begin{equation}
\E[S^{\frac{2}{\beta}}]<\infty\,.\label{e.MomentAssumption}%
\end{equation}
Then the process of propagation losses $\mathcal{N}$  experienced by the typical user
 is a non-homogeneous Poisson point
process on $\mathbb{R}^{+}$\ with intensity measure
\begin{equation}
\Lambda\left(  \left[  0,t\right)  \right)  :=\E\left[  \mathcal{N}
\left(  \left[  0,t\right)  \right)  \right]  =at^{\frac{2}{\beta}}
\label{e.Displacement}%
\end{equation}
where
\begin{equation}
a:=\frac{\lambda\pi \E[S^{\frac{2}{\beta}}]}{K^{2}} \label{e.ModifiedIntensity}%
\end{equation}
\end{lemma}
\begin{IEEEproof} 
The point process $\mathcal{N}$  may be viewed as a transformation of the point
process~$\Phi$ by the probability
kernel
\[
p(x,A)=\Prob\left(  \frac{l\left(  x\right)  }{S}\in A\right)  ,\quad
x\in\mathbb{R}^{2},A\in\mathcal{B}\left(  \mathbb{R}^{+}\right)
\]
By the displacement theorem~\cite[Theorem 1.10]{FnT},
the point process $\mathcal{N}$\ is Poisson on $\mathbb{R}^{+}$\ with
intensity measure
\begin{align*}
\Lambda\left(  \left[  0,t\right)  \right)   
&  =\lambda\int_{\mathbb{R}^{2}}\Prob\left(  \frac{l\left(  x\right)  }{S}%
\in\left[  0,t\right)  \right)  dx\\
&  =\lambda\int_{\mathbb{R}^{2}\times\mathbb{R}^{+}}\Ind\left\{  \frac{l\left(
x\right)  }{s}<t\right\}  dx\Prob_{S}\left(  ds\right) \\
&  =\lambda\int_{\mathbb{R}^{+}}\frac{\pi\left(  st\right)  ^{\frac{2}{\beta}%
}}{K^{2}}\Prob_{S}\left(  ds\right)  =\frac{\lambda\pi \E\left[  S^{\frac{2}{\beta
}}\right]  }{K^{2}}t^{\frac{2}{\beta}}%
\end{align*}
which completes the proof.
\end{IEEEproof}

\begin{remark}
Note that the distribution of $\calN$ is {\em invariant with respect
  to the distribution of the shadowing/fading $S$ having  same given
  value of the   moment~$\E[S^{2/\beta}]$}.
This means that the infinite Poisson network with an arbitrary
  shadowing $S$   
is perceived at a given location statistically in the same manner
 as an ``equivalent'' infinite Poisson  with ``constant shadowing'' 
 equal to $s_{const}=(\E[S^{2/\beta}])^{\beta/2}$ (to have the same moment of
 order $2/\beta$). The model with such a ``constant shadowing'' 
boils down to the model without shadowing ($S\equiv 1$) and the constant
$K$ replaced by $\tilde K=K/\sqrt{\E[S^{2/\beta}]}$.
\end{remark}
The above result requires condition (\ref{e.MomentAssumption}), which
is satisfied by the usual models, as e.g.
the log-normal shadowing or Rayleigh fading, 
and we tacitly assume it throughout the paper.

\section{Convergence results under log-normal shadowing}
\label{s.convergence}
In this section we derive a powerful convergence result rigorously showing that
the infinite Poisson model can be used to analyse the characteristics
of the typical user in the context of any {\em fixed (deterministic!)
placement of base stations}, meeting some \emph{empirical homogeneity}
condition,  provided there is {\em sufficiently strong log-normal shadowing}. 
\subsection{Model description}
Let $\phi=\{X_i\}_{i\in \mathbb{N}}$ be a locally finite deterministic
point pattern  on $\R^2$ 
and $B_0(r)$ the ball of radius $r$,
centered at the origin.
For $0 < \lambda < \infty$, as $r\rightarrow \infty$ we require the
{\em empirical homogeneity condition}
\begin{equation}\label{lambdaDef}
 \frac{\phi(B_0(r))}{\pi r^2} \rightarrow \lambda.
\end{equation}
Note that the above condition is satisfied by any lattice
(e.g. hexagonal) pattern~$\phi$, as well as by almost any realization 
of an arbitrary ergodic point process.

Let the  shadowing  $S_i^{(\sigma)}$ between the station $X_i$ and the
origin  be iid (across $i$) log-normal random variables 
\begin{equation}\label{e.lognormalRV}
S_i^{(\sigma)}=\exp(-\sigma^2/2+\sigma Z_i),
\end{equation}
where $Z_i$ are standard normal random variables.
Note that for such $S_i^{(\sigma)}=S^{(\sigma)}$, we have
 $\E[S]=1$ and $\E\left[
  (S^{(\sigma)})^{2/\beta}\right]=\exp[\sigma^2(2-\beta)/\beta^{2}]$.~\footnote{%
\label{ftnte}
This also means that the path-loss from a given station $X$
 expressed in dB, i.e., $\dB(L_X(y))$, where
$\dB\left(  x\right):=10\times\log_{10}\left(x\right)\dB$, is a Gaussian
random variable with  standard deviation $\logsd=\sigma10/\log10$,
called {\em logarithmic standard deviation of the shadowing}.}

Consider the distance-loss model~(\ref{e.DistanceLoss}) with the
constant $K$ replaced by the function of $\sigma$
\begin{equation}
K(\sigma)=K \exp\left(-\frac{\sigma^2(\beta-2)}{2\beta^2}\right),
\end{equation}
where  $K>0$ and $\beta>2$.

As in Section~\ref{PoissonModel1}, we consider the point process on $\mathbb{R}^+$
of  propagation
  losses  experienced by the typical user with respect to the stations
in~$\phi$
\begin{equation}
 \calN^{(\sigma)}:= \left\{\frac{ K(\sigma)^{\beta} |X_i|^{\beta}
 }{S_i^{(\sigma)}}: X_i\in\phi \right\}\,.
\end{equation}
We consider also the analogous process of  propagation
  losses  
\begin{equation}
 \bar\calN^{(\sigma)}:= \left\{\frac{ K(\sigma)^{\beta} |X_i|^{\beta}
 }{S_i^{(\sigma)}}: a_\sigma < |X_i| <b_\sigma, X_i\in\phi \right\}\,,
\end{equation}
where the stations in $\phi$ that are closer than $a_\sigma$ and farther
than $b_\sigma$ are ignored, for all sequences
$0\le a_\sigma< b_\sigma \le \infty$ satisfying
\begin{equation}\label{a_n}
 \frac{\log( a_\sigma)}{\sigma^2} \rightarrow 0,
\end{equation}
\begin{equation}\label{b_n}
 \frac{\log(b_\sigma)}{\sigma^2} \rightarrow \infty.
\end{equation}

\subsection{Main result}
We present now our main convergence result.
\begin{theorem}\label{mainResult}
Given homogeneity condition (\ref{lambdaDef}), then  
$\calN^{(\sigma)}$
converges weakly as $\sigma\rightarrow \infty$ to the
 Poisson point process on $\R^+$ with the intensity measure $\Lambda$
 given by~(\ref{e.Displacement}) with $a=\lambda\pi/K^{2}$. 
Moreover, $\bar\calN^{(\sigma)}$
also converges weakly $(\sigma\rightarrow \infty)$ 
to the Poisson point process with the same intensity measure,
provided conditions~(\ref{a_n}) and~(\ref{b_n}) are satisfied.
\end{theorem}



The  proof of Theorem~\ref{mainResult} is deferred to Appendix~\ref{App1}.

\begin{remark}\label{remark2}
The above result, in conjunction with Lemma~\ref{l.ModifiedIntensity},
 says that {\em infinite Poisson model can be used
to approximate the characteristics of the typical user for a very general
class of homogeneous pattens of base stations}, including the standard
hexagonal one. The second statement of this result says that this
approximation  remains valid for {\em sufficiently large but finite patterns}.
\end{remark}

\begin{remark}
The distance-loss model~(\ref{e.DistanceLoss}) suffers
from having a singularity at the origin. This issue is often
circumvented by some  appropriate modification of the distance-loss function 
within a certain distance from the origin.  
The  second statement of Theorem~\ref{mainResult} with
$a_\sigma=const>0$  shows that such a modification is not significant
in the Poisson approximation.
\end{remark}

\begin{figure}[t!]
\begin{center}
\centerline{\includegraphics[width=0.8\linewidth]{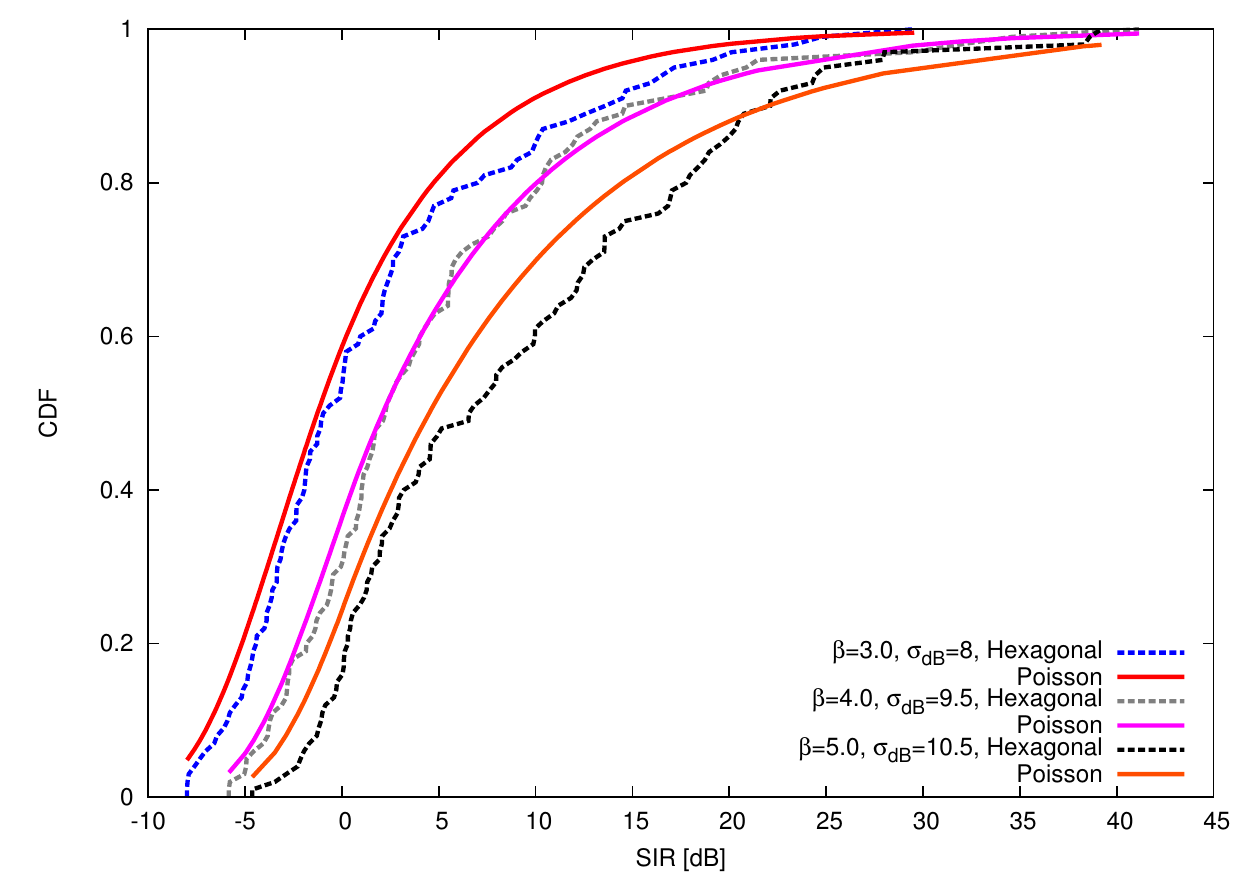}}
\caption{Empirical CDF of SIR simulated in hexagonal network with shadowing
and their Poisson approximations.
\label{f.SIR}}
\end{center}
\end{figure}

To illustrate Theorem~\ref{mainResult} and obtain some insight into
the speed of convergence we used  Kolmogorov-Smirnov (K-S) test,
cf~\cite{Williams2001}, to compare the cumulative distribution
functions (CDF) of the SIR of the typical user (see
Section~\ref{sss.SIR}) in the infinite 
Poisson model versus hexagonal one consisting of   $30\times30=900$
stations  on a torus. We found that for 
9/10 realizations of the network shadowing
the K-S test does not allow to distinguish
the empirical (obtained from simulations)  CDF of the SIR 
from the CDF of SIR evaluated in the infinite  Poisson model 
with the  critical $p$-value fixed to $\alpha=10\%$
provided  $\logsd\gtrsim 10\text{dB}$.
On Figure~\ref{f.SIR} we present a few examples of these CDF's.
Similar numerical study was done for the CDF of the path-loss with respect
to the serving station in~\cite{propag_extabst}.

\section{Poisson model revisited}\label{PoissonModel2}
We now return to the model outlined in Section \ref{PoissonModel1} and
investigate it further.  We show that Lemma~\ref{l.ModifiedIntensity} is very useful in
studying random quantities of the network, independently of the
shadowing/fading distribution. We do not pretend to develop a complete
theory, which is well beyond the scope of this paper. As a general
remark let us emphasize only that many, already known, results (cf
{\em Related work}  in the Introduction), 
were originally derived under specific
assumptions on the distribution of shadowing and/or
fading. By Lemma~\ref{l.ModifiedIntensity} they necessarily remain valid for a
general distribution, with the appropriate specification of the value of the
moment~$\E[S^{2/\beta}]$.

Our specific goal in this section is to present a new representation of
distribution function of the SINR, which in consequence will
allow us to study the spectral and energy efficiency --- two
important engineering characteristics of the network. To the best of
our knowledge, they have not yet been studied in the Poisson model.
Moreover, as another illustration of our convergence result, we will
compare these characteristics in Poisson and hexagonal network, with this
latter  studied by  simulations.
 
Unless otherwise specified all results of this section regard Poisson
model of Section~\ref{PoissonModel1}.

\subsection{Path-loss and SIR}
The novelty of our approach consists in representing the SINR in terms
of the path-loss to the serving station and the respective SIR,
rather than, as usual, the interference.

\subsubsection{Path loss} The weakest propagation loss denoted by
\begin{equation}
L=\inf_{i\in\mathbb{N}}L_{X_{i}} \label{e.MinLoss}%
\end{equation}
is often called \emph{path-loss factor}. 
It may be interpreted
as the propagation loss with the \emph{serving} base station (that is, the one
with strongest received power). Note by~(\ref{e.Displacement})
that the number of points of the process $\mathcal{N}=\left\{  L_{X_{i}%
}\right\}  _{i\in\mathbb{N}}$\ is almost surely finite in any finite
interval. Hence, the above infimum is almost surely achieved for some base
station; that is $L=\min_{i\in\mathbb{N}}L_{X_{i}}$.
Moreover, Lemma~\ref{l.ModifiedIntensity} allows us to conclude the
following characterisation of the distribution of~$L$, which follows
immediately from the well known expression for void probabilities of
Poisson process $\P\{L\ge t\}=\P\{\calN([0,t))=0\}=\exp[-\Lambda([0,t))]$.
\begin{corollary}
\label{c.MinCDF}
The CDF of $L$ is equal to 
$\P\{\, L \le  t\,\}  =1-\exp[-at^{2/\beta}]$
where $a$ is given by~(\ref{e.ModifiedIntensity}).
Consequently, the probability density of $L$ is 
given by
\begin{equation}
\Prob_{L}\left(  ds\right)  =\frac{2a}{\beta}t^{\frac{2}{\beta}-1}e^{-at^{\frac
{2}{\beta}}}dt\,.\label{e.LossPDF}%
\end{equation}
This corresponds to a \emph{Fr\'{e}chet distribution} with shape parameter $\frac{2}{\beta}$\ and
scale parameter $a^{\frac{\beta}{2}}$.
\end{corollary}

\subsubsection{Interference factor and SIR}
\label{sss.SIR}
The \emph{interference factor} is defined by
\[
f=L\sum_{i\in\mathbb{N}}\frac{1}{L_{X_{i}}}-1
\]
where $L$\ is the path-loss factor. It may be interpreted as the interference
to signal ratio; that is the \emph{inverse of the SIR}.
Introducing  $I=\sum_{i\in\mathbb{N}}\frac{1}{L_{X_{i}}}$
(can be interpreted as the total received power) 
we can write $f=L I-1$. Also, $f =L I^{\prime}$
where
\[
I^{\prime}:=I-\frac{1}{L}=\sum_{i\in\mathbb{N}}\frac{1}{L_{X_{i}}}-\frac{1}{L}%
\]
 may be viewed as the interference.

Define
\begin{equation}
\varphi_{\beta}\left(  z\right)  :=e^{-z}+z^{\frac{2}{\beta}}\gamma\left(
1-\frac{2}{\beta},z\right)  \label{e.LaplaceF0}%
\end{equation}
where $\gamma(\alpha,z)=\int_{0}^{z}t^{\alpha-1}e^{-t}dt$\ is\ the
lower\ incomplete\ gamma\ function.
Another representation of the function $\varphi_{\beta}(z)$, 
used in evaluation of the Laplace transform of the $f$
(cf proof of Corollary~\ref{p.LaplaceFandJoint}) is given in Appendix~\ref{App2}.

Here is a key technical result for our approach.
\begin{proposition}
\label{p.LaplaceF_CondL} The Laplace transform of the interference factor
conditional to the path-loss factor is equal to 
\begin{equation}
\E\left[  e^{-zf}|L=s\right]  =e^{-a\left[  \varphi_{\beta}(z)-1\right]
s^{\frac{2}{\beta}}}\,. \label{e.LaplaceF_CondL}%
\end{equation}
\end{proposition}

\begin{IEEEproof}
We have
\begin{equation}
\E\left[  e^{-zf}|L=s\right]  =\E[e^{-zsI^{\prime}}|L=s] \label{e.LaplaceF1}\,.
\end{equation}
Observe that  $I^{\prime}$ 
is a shot-noise  associated to the point process
obtained  from the Poisson point process $\calN$ of the propagation
losses by suppressing its smallest (nearest to the origin) point,
which is located precisely at $L$. By the well-known property
of the Poisson process, given $L=s$, this new process is also Poisson
on $(s,\infty)$ with intensity measure
\[
{M}\left(  \left(  s,t\right)  \right)  =\Lambda\left(  \left[
0,t\right)  \right)  -\Lambda\left(  \left[  0,s\right]  \right)
,\quad t\in\left(  s,+\infty\right)\,.
\]
Thus
\begin{equation}
{M}\left(  dt\right)  =\Lambda\left(  dt\right)  =\frac{2a}{\beta
}t^{\frac{2}{\beta}-1}dt,\quad\text{on }\left(  s,+\infty\right) 
\label{e.DisplacementConditioned}%
\end{equation}
and by~\cite[Proposition 1.5]{FnT}
\begin{align}
\E\left[  e^{-zsI^{\prime}}|L=s\right]   &  =\exp\left[  \int_{s}^{\infty
}(e^{-\frac{zs}{t}}-1){M}\left(  dt\right)  \right] \nonumber\\
&  =\exp\left[  -\frac{2a}{\beta}\int_{s}^{\infty}(1-e^{-\frac{zs}{t}%
})t^{\frac{2}{\beta}-1}dt\right] \,. \label{e.LaplaceF2}%
\end{align}
Now  we calculate the integral on the right-hand side of the above
equation by making the  change of variable $u:=\frac{zs}{t}$, that is
\begin{eqnarray*}
\lefteqn{\int_{s}^{\infty}(1-e^{-\frac{zs}{t}})t^{\frac{2}{\beta}-1}dt}\\
&=&  \int_{z}^{0}(1-e^{-u})\left(  \frac{zs}{u}\right)  ^{\frac{2}{\beta}%
-1}\frac{-zs}{u^{2}}du\\
&=&  \left(  zs\right)  ^{\frac{2}{\beta}}\int_{0}^{z}(1-e^{-u})u^{-\frac
{2}{\beta}-1}du\\
&=&  s^{\frac{2}{\beta}}\frac{\beta}{2}\left[  -1+\varphi_{\beta}\left(
z\right)  \right]
\end{eqnarray*}
where the third equality is obtained by integration by parts. Combining the above equation with Equations~(\ref{e.LaplaceF1}),
(\ref{e.LaplaceF2}) gives the final result.
\end{IEEEproof}

The following results can be derived from
Proposition~\ref{p.LaplaceF_CondL}, although we will not use them in
the remaining part of the paper.
\begin{corollary}
\label{p.LaplaceFandJoint}
The joint distribution of $L$\ and $f$ is
characterized by
\[
\E\left[  \Ind\left\{  L\geq u\right\}  e^{-zf}\right]  =\frac{1}{\varphi_{\beta
}\left(  z\right)  }e^{-au^{\frac{2}{\beta}}\varphi_{\beta}\left(  z\right)
},\quad z\in\mathbb{R}^{+}\,.%
\]
The unconditional Laplace transform of the interference factor is
\begin{equation}
\E\left[  e^{-zf}\right]  =\frac{1}{\varphi_{\beta}\left(  z\right)  },\quad
z\in\mathbb{R}^{+}\,. \label{e.LaplaceF}%
\end{equation}
\end{corollary}
\begin{IEEEproof}
For the first statement we have
\begin{align*}
\E\left[  \Ind\left\{  L\geq u\right\}  e^{-zf}\right]   &  =\int_{u}^{\infty
}\E\left[  e^{-zf}|L=s\right]  \Prob_{L}\left(  ds\right) \\
&  =\int_{u}^{\infty}e^{-a\left[  \varphi_{\beta}(z)-1\right]  s^{\frac
{2}{\beta}}}\frac{2a}{\beta}s^{\frac{2}{\beta}-1}e^{-as^{\frac{2}{\beta}}}ds\\
&  =\frac{1}{\varphi_{\beta}\left(  z\right)  }e^{-au^{\frac{2}{\beta}}%
\varphi_{\beta}\left(  z\right)  }%
\end{align*}
where for the second equality we use~(\ref{e.LossPDF})
and~(\ref{e.LaplaceF_CondL}). This completes the proof of the first statement.
For the second statement, 
conditioning on the value of the path-loss factor yields
\[
\E\left[  e^{-zf}\right]  =\int_{\mathbb{R}^{+}}\E\left[  e^{-zf}|L=s\right]
\Prob_{L}\left(  ds\right)
\]
where $\Prob_{L}\left(  \cdot\right)  $\ is\ the\ distribution\ of\ $L$ given
by~(\ref{e.LossPDF}).\ Using the above equation and~(\ref{e.LaplaceF_CondL}),
we ascertain
\begin{align*}
\E[e^{-zf}]  &  =\int_{0}^{\infty}e^{-a\left[  \varphi_{\beta}(z)-1\right]
s^{\frac{2}{\beta}}}\frac{2a}{\beta}s^{\frac{2}{\beta}-1}e^{-as^{\frac
{2}{\beta}}}ds\\
&  =\frac{2a}{\beta}\int_{0}^{\infty}s^{\frac{2}{\beta}-1}e^{-a\varphi_{\beta
}(z)s^{\frac{2}{\beta}}}ds\\
&  =\frac{1}{\varphi_{\beta}(z)}\int_{0}^{\infty}a\varphi_{\beta}(z)\frac
{2}{\beta}s^{\frac{2}{\beta}-1}e^{-a\varphi_{\beta}(z)s^{\frac{2}{\beta}}}ds\\
&  =
\frac{1}{\varphi_{\beta}(z)}%
\end{align*}
where the third equality stems from $\varphi_{\beta}(z)\neq
0$\, which follows from~(\ref{e.PhiPositive}) in Appendix~\ref{App2} 
and the assumption $\beta>2$.
\end{IEEEproof}

\begin{remark}
\label{r.SIRexplicit}
For $t\ge 1$ the (complementary) 
CDF of the SIR admits the following explicit expression
\begin{equation}\label{e.SIRexplicit}
\P\{\,\text{SIR}\ge t\,\}=\P\{f\le 
1/t\}=\frac{t^{-2/\beta}}{C'(\beta)}\,
\end{equation}
where 
$$C'(\beta)=\frac{2\pi}{\beta\sin(2\pi/\beta)}
=\frac{2\Gamma(2/\beta)\Gamma(1-2/\beta)}{\beta}\,$$
and  $\Gamma(z)=\int_0^\infty e^{-z}z^{t-1}dt$ 
it the complete Gamma function.
It was proved in~\cite{Dhillon2012} assuming exponential  distribution of
$S$. But in the infinite Poisson model, SIR is   invariant
with respect to $\lambda$, and $K$ and  consequently by Lemma~\ref{l.ModifiedIntensity}
it is also invariant with respect to the value of  $\E[S^{2/\beta}]$.
Hence the result remains valid for {\em arbitrary distribution} of
$S$, provided $\E[S^{2/\beta}]<\infty$. Other results
of~\cite{Dhillon2012}, including  these involving 
the superposition of independent Poisson models (called $K$-tier
cellular network model) can also be appropriately generalized using Lemma~\ref{l.ModifiedIntensity}.
\end{remark}

\subsection{Distribution of SINR, spectral and energetic efficiency}
\subsubsection{SINR}
In this section we are primarily interested in the {\em signal to interference and noise ratio}
\begin{equation}
\mathrm{SINR}=\frac{\frac{1}{L}}{N+\left(  \sum_{i\in\mathbb{N}}\frac
{1}{L_{X_{i}}}-\frac{1}{L}\right)  }=\frac{1}{NL+f}, \label{e.SINR}%
\end{equation}
where $N$ is the noise power. We will show how  CDF of the SINR can be
evaluated, which opens a way for the study of functionals of SINR.

\subsubsection{Spectral efficiency} 
An important characteristic of 
a wireless cellular network is its \emph{spectral efficiency},
defined in the simplest case of additive white Gaussian noise (AWGN)
and the optimal  theoretical link performance 
as  
$$\calS:=\log\left(  1+\mathrm{SINR}\right)\,.$$
It tells us how many  bits per second and per Hertz  can be sent to the typical user of the
network. 

\subsubsection{Energy efficiency}
Up to now we have considered a unit transmitted power.
Assume now that base stations transmit some power $P\ge0$. In fact,
in order for a base station to be able to transmit this power to the
mobile, it needs to be powered (i.e., consumes energy per second)
at the level $P'>P$. Following~\cite{Richter2009},
assume that these two quantities are   related  through a simple linear relation 
$P' = cP + d$ for some positive constants $c$ and $d$. 
The  {\em energy efficiency} is
defined by 
$$\calE:=\calE (P) = \frac{W\,\log\left(1+\frac{1}{NL/P+f}\right)}{cP+d}\,.$$
where $W$ is the bandwidth expressed in Hz. Thus $\calE$ is equal
to $W\log(1+\text{SINR}(P))/P'$, where $\text{SINR}(P)$ takes into account the
transmitted power, and $P'$ is the consumed power.
It tells us how many bits per second per Watt of consumed power
can be sent by a base station to the typical user. Note that
 $\calE(0)=\calE(\infty)=0$ and thus $\calE(P)$ admits a
non-trivial optimization in~$P$.

\subsubsection{Evaluation of the CDF of  the SINR}
\label{sss.SINR}
Proposition~\ref{p.LaplaceF_CondL} in conjunction with
Corollary~\ref{c.MinCDF} completely characterizes the
joint distribution of $L$\ and $f$. It does not, however, allow for an
explicit expression for the CDF of the SINR in the whole domain
(cf Remark~\ref{r.SIRexplicit}).
In this section we  describe a practical way for numerical computation
of this CDF. We will use it to study the spectral and energy efficiency.
In fact, we will compute the CDF of the random variable $Y:=NL+f$,
which is sufficient, in view of~(\ref{e.SINR}).

\begin{proposition}
\label{p.INSR_CDF}The cumulative distribution function of $Y$ is given by
\begin{equation}
\Prob\left(  Y<x\right)  =\int_{0}^{\infty}F_s(x-Ns)\Prob_{L}\left(  ds\right)  \label{e.INSR_CDF}%
\end{equation}
where $\Prob_{L}\left(  ds\right)  $\ is given by~(\ref{e.LossPDF}), and
$F_{s}\left(  y\right)  :=\Prob\left(  f<y|L=s\right)$ may be
expressed by
\begin{equation}\label{eq:Bromwich}
F_s(y)=1-\frac{2e^{\gamma t}}{\pi}
\int_{0}^\infty \calR\left(\overline\calL_{F_s}
(\gamma+iu)\right)\cos ut\, d u\,,
\end{equation}
where $\gamma>0$ is an arbitrary constant, $\calR$ is the real part if it
complex argument 
 and
\begin{equation}\label{e.LTFs}
\overline{\mathcal{L}}_{F_{s}}\left(  z\right)  =\frac1z-\frac{1}{z}e^{-a\left[  \varphi_{\beta
}(z)-1\right]  s^{\frac{2}{\beta}}}\,%
\end{equation}
with $\varphi_{\beta}\left(  \cdot\right)  $\ is given by~(\ref{e.LaplaceF0}).
\end{proposition}

\begin{IEEEproof}
The expression~(\ref{e.INSR_CDF}) is trivial.
The expression~(\ref{eq:Bromwich}) of  the (conditional) CDF $F_{s}$ 
is based on the Bromwich contour inversion integral of the Laplace
transform $\overline{\mathcal{L}}_{F_{s}}$ of $1-F_{s}$. The
expression~(\ref{e.LTFs}) for this latter follows
from~(\ref{e.LaplaceF_CondL}).
\end{IEEEproof}

\begin{remark}
As shown in~\cite{AbateWhitt95}, the integral in~(\ref{eq:Bromwich}) 
can be numerically evaluated using the trapezoidal rule, with the
parameter~$\gamma$ allowing control of the approximation error.
The function  $F_s$  may be also
retrieved from~(\ref{e.LaplaceF_CondL})  using other inversion techniques.
\end{remark}

\begin{remark}
\label{r.SINRexplicit}
For $t\ge 1$ the (complementary) 
CDF of the SINR admits the following  expression
\begin{equation}\label{e.SINRexplici}
\P\{\text{SINR}\ge t\}=\frac{2t^{-2/\beta}}{\Gamma(1+\frac{2}{\beta})}
\int_0^\infty \!\!\!\!\!\!re^{-r^2\Gamma(1-2/\beta)- Na^{-\beta/2} r^\beta}dr
\end{equation}
with $a$ given by~(\ref{e.ModifiedIntensity}).
It follows from \cite[Theorem~1]{Dhillon2012} (where 
the exponential distribution is assumed for 
$S$) and our  Lemma~\ref{l.ModifiedIntensity}.
Note on Figure~\ref{f.SINR} that the expression is not valid for $t<1$.
\end{remark}

\subsubsection{Numerical examples}
\label{sss.Numerical}
In what follows we assume the following parameter values.
If not otherwise specified we take log-normal shadowing
with logarithmic standard deviation $\logsd=12$dB
(cf the footnote on page~\pageref{ftnte}).
The parameters of the distance-loss model are  $K=4250$Km$^{-1}$,  $\beta=3.52$
(which corresponds to the  COST-Hata model for urban environment with the BS height 30m, 
mobile height 1.5m, carrier frequency 805MHz and penetration loss 19dB). 
For the hexagonal network model simulations, we consider  900 base stations (30x30) on a torus,
with the cell radius  $R=0.26$Km (i.e., the surface of the hexagonal
cell is equal to this of the disk  of radius $R$). 
System bandwidth is $W=10$MHz, noise power $N=-93$dBm.

\begin{figure}[t!]
\begin{center}
\centerline{\includegraphics[width=0.8\linewidth]{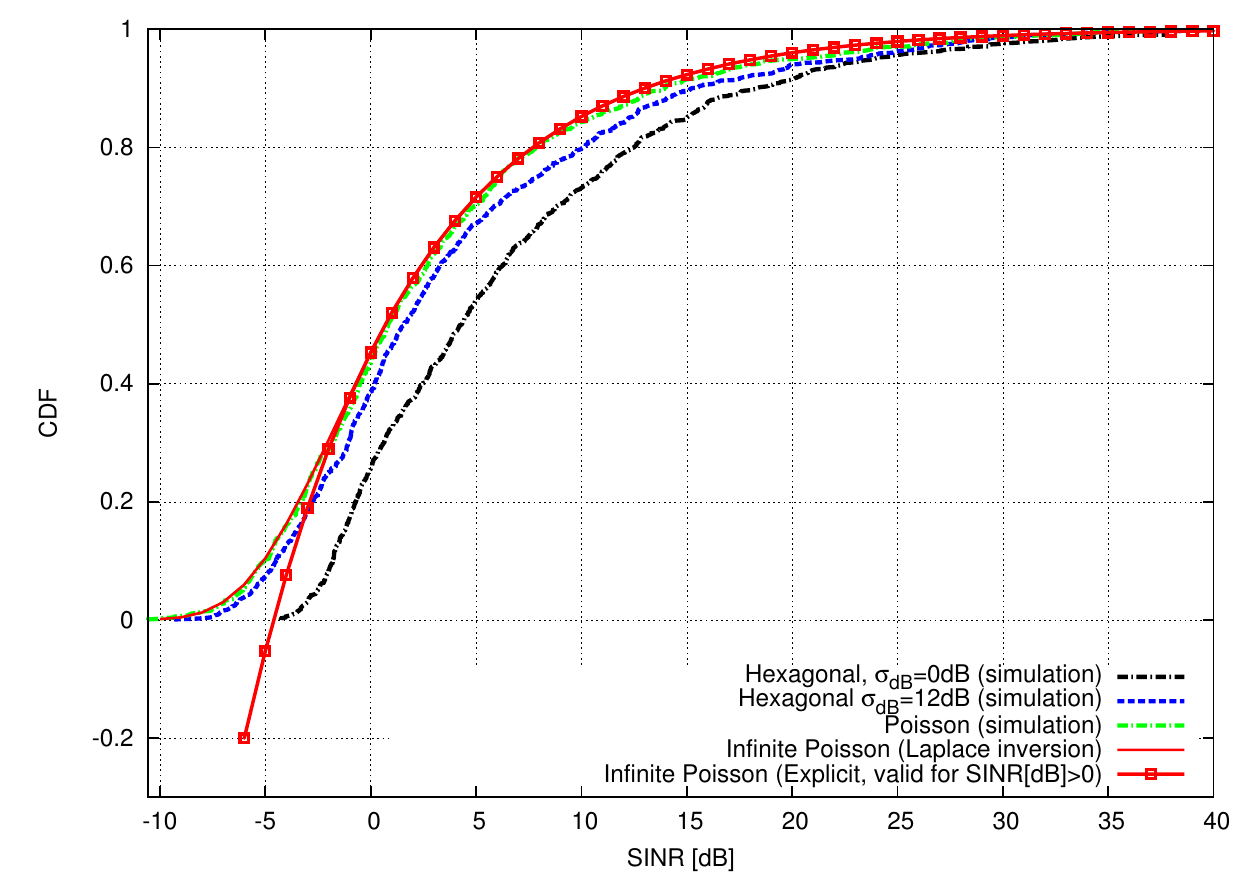}}
\caption{CDF of the SINR for (finite) hexagonal network 
 with and without  ($\logsd=0\text{dB}$) shadowing, 
finite Poisson as well as for the
 infinite Poisson model. The last in the legend curve corresponds to~(\ref{e.SIRexplicit}).
\label{f.SINR}}
\end{center}
\end{figure}

Figure~\ref{f.SINR} shows the CDF of the SINR assuming the base
station  power $P=58.5$dBm. We obtain results by simulations for finite hexagonal network on the torus with and
without shadowing as well as the finite Poisson network on the same
torus. For the infinite Poisson 
model we use our method based on the inversion of the Laplace
transform. For comparison, we plot also the curve corresponding the
explicit expression~(\ref{e.SIRexplicit}), 
which is valid only for $\text{SINR}>1$.

\begin{figure}[t!]
\begin{center}
\centerline{\includegraphics[width=0.8\linewidth]{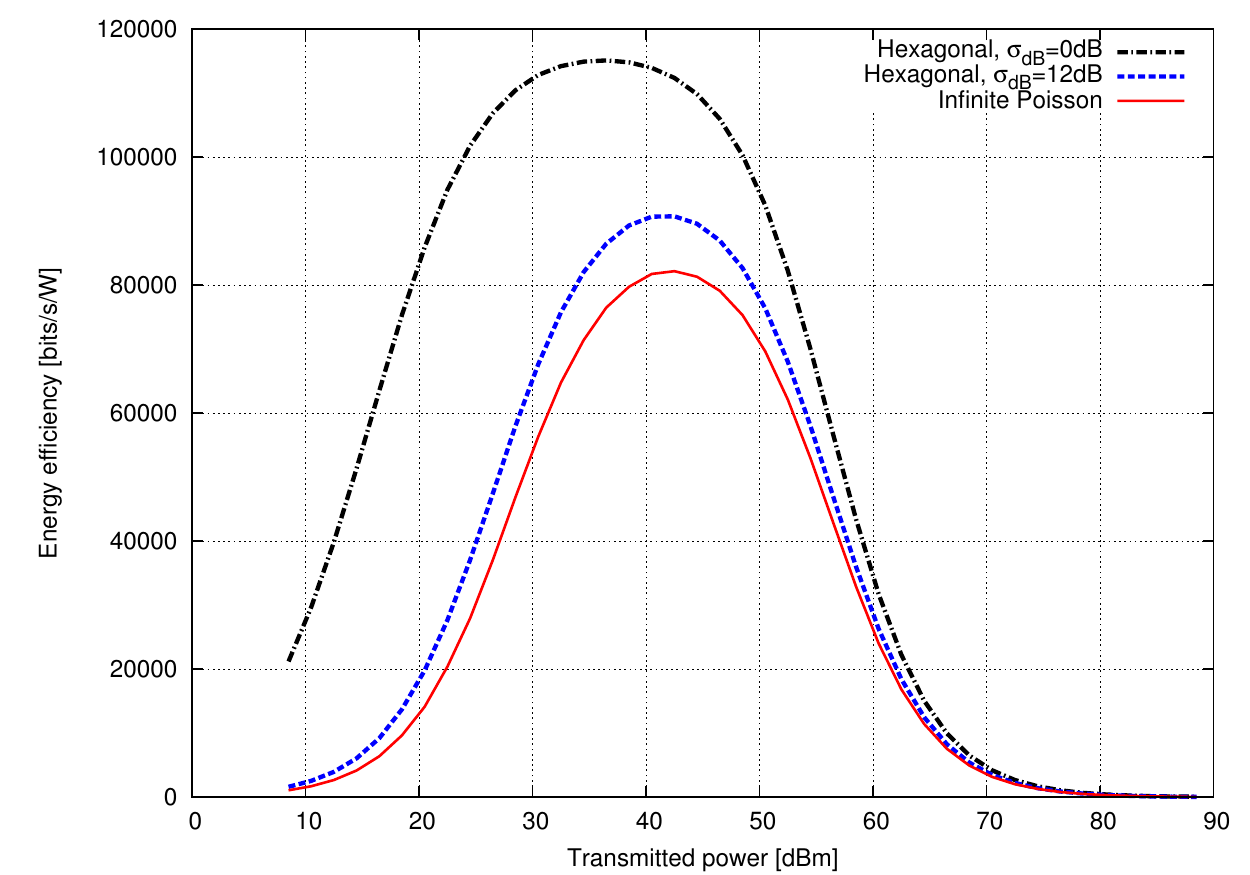}}
\caption{Mean energy efficiency of finite hexagonal model
with and without ($\logsd=0\text{dB}$) shadowing as well as the infinite Poisson
in function of the transmitted power.
\label{f.E}}
\end{center}
\end{figure}

Figure~\ref{f.E} shows the expected energy efficiency $\E(\calE)$ of the finite hexagonal model
with and without shadowing as well as the infinite Poisson,
assuming the affine relation between consumed and emitted power
with constants $c =21,45$ and $d = 354.44$W.

On both figures we see that the infinite Poisson model gives a reasonable
approximation of the (finite) hexagonal network provided the
shadowing is high enough. In particular, the value of the transmitted
power at which the hexagonal network with the shadowing attains the
maximal expected energy efficiency is very well predicted by the Poisson model.

\section*{Conclusions}
We present a powerful mathematical convergence result 
rigorously justifying the fact
that {\em any} actual (including regular hexagonal) network 
is perceived by a typical user as an equivalent (infinite) Poisson network, 
provided log-normal shadowing variance is sufficiently high.
Good approximations are
obtained for logarithmic standard deviation of the shadowing  greater
than approximately 10dB, which is a  realistic assumption in many urban scenarios.
Moreover, the  equivalent infinite Poisson representation is invariant
with respect to an additional fading distribution.
Using this representation we study the distribution
of the SINR of the typical user.
In particular, we evaluate and optimize 
the  mean energy efficiency
as a function of the base station transmit power.

\appendix

\subsection{Proof of Theorem~\ref{mainResult}}
\label{App1}
In order to simplify the notation we set $n:=\sigma^2$. Moreover,
without loss of generality we assume that $n$ takes positive integer
values. Also, it is more convenient to study the  point
process of propagation losses on the logarithmic scale, which we do in
what follows. In this regard denote by $\Lambda_{\log}$  the image of the
measure $\Lambda$ given by~(\ref{e.Displacement}) with
$a=\lambda\pi/K^2$ through the
logarithmic mapping 
$$\Lambda_{\log}((-\infty,s]):= \int_{R^+}1(\log(t)\le s)\,\Lambda(dt)
=\frac{\pi}{K^2} e^{\frac{2s}{\beta}}$$
for $s\in\R$.

For a given $i$, we first observe by~(\ref{e.lognormalRV})  that
\begin{align}\nonumber
\nu_n (s,|X_i|)&:=\Prob \left[  \log\left(\frac{K(n)^{\beta} |X_i|^{\beta} }{S_i^{(n)}} \right) \leq
s\right ] \\
&=\Prob \left[Z_i \leq 
\frac{s-\beta\log(K|X_i|)-n/\beta}{\sqrt{n}} \right ]\nonumber \\
&=G \left[\frac{s-\beta\log(K|X_i|)-n/\beta}{\sqrt{n}} \right ],\label{v_nFinal}
\end{align}
where $G$ is the CDF of the standard Gaussian random variable.

Let $B_0(r)=\{x\in\R^2:|x|<r\}$. 
We now need to derive two results.
\begin{lemma}\label{Lemma1}
\begin{align*}
\lim_{n\rightarrow \infty} \int_{B_0(b_n)\setminus B_0(a_n)} \nu_n(s,|x|)dx &=\lim_{n\rightarrow \infty}\int_{\R^2}\nu_n(s,|x|)dx \\
& =\Lambda_{\log}((-\infty,s])
\end{align*}
provided that $a_n$ and $b_n$ satisfy (\ref{a_n}) and (\ref{b_n}).
\end{lemma}
\begin{IEEEproof}
We start with the integral over $\R^2$ in polar form
\begin{align*}
 \int_{\R^2}\nu_n (s,|x|)dx  
&= 2\pi
\int_0^{\infty}rG \left[\frac{s-\beta\log(Kr)-n/\beta}{\sqrt{n}} \right
]dr
\\
\end{align*}
Introduce a change of variables
\begin{align*}
 t&=\frac{s-\beta\log(Kr)-n/\beta}{\sqrt{n}}, \qquad \qquad
&dt=-\frac{\beta}{\sqrt{n}}\frac{dr}{r},\\
 r&=\frac{1}{K}\exp\left[ \frac{1}{\beta}\left(s-t\sqrt{n} -n/\beta \right)
\right],  &dr=-\frac{\sqrt{n}r}{\beta}dt,
\end{align*}
hence
\begin{align*}
&2\pi
\int_0^{\infty}rG \left[\frac{s-\beta\log(Kr)-n/\beta}{\sqrt{n}}
\right ]dr \\
&=2\pi \int_{-\infty}^{\infty} \frac{1}{K^2}\exp\left[
\frac{2}{\beta}\left(s-t\sqrt{n} -n/\beta \right) \right] 
G (t)\frac{\sqrt{n}}{\beta}dt \\
&=2\pi \frac{\sqrt{n}}{\beta} \frac{\exp\left[ \frac{2}{\beta}\left(s
-n/\beta \right) \right]}{K^2}
 \int_{-\infty}^{\infty} \exp\left[ \frac{-2t\sqrt{n}}{\beta} \right] 
G (t)dt 
\end{align*}
Noting that $G (-t)=1-G (t):=\bar{G }(t)$, apply integration by parts on
the above integral
\begin{align*}
 \int_{-\infty}^{\infty} \exp\left[ \frac{-2t\sqrt{n}}{\beta} \right] 
G (t)dt &=
 \int_{-\infty}^{\infty} \exp\left[ \frac{2t\sqrt{n}}{\beta} \right] 
\bar{G }(t)dt
\\ &=\left.\frac{\beta}{2\sqrt{n}} \exp\left[ \frac{2t\sqrt{n}}{\beta}
\right]
 \bar{G }(t)\right| _{-\infty}^{\infty} \\
+&\frac{\beta}{2\sqrt{n}}\int_{-\infty}^{\infty} \exp\left[
\frac{2t\sqrt{n}}{\beta} \right]  g (t)dt,
\end{align*}
where $g(t)$ is the normal probability density, and, via the inequality 
\begin{equation}\label{errIneq}
\int_y^{\infty}e^{-x^2}dx<\frac{e^{-y^2}}{(1+y)}, \quad y\geq 0\,,
\end{equation}
(cf \cite[Section~7.8]{DLMF:2010})
the final integrated term on the left vanishes. Hence
\begin{align}
 &\int_{-\infty}^{\infty} \exp\left[ \frac{-2t\sqrt{n}}{\beta} \right] 
G (t)dt \nonumber\\
&=
\frac{\beta}{2\sqrt{n}}\int_{-\infty}^{\infty} \exp\left[-\frac{t^2}{2}
+\frac{2t\sqrt{n}}{\beta} \right] \frac{dt}{\sqrt{2\pi}}\nonumber\\
&=\frac{\beta}{2\sqrt{n}} \exp\left[\frac{2n}{\beta^2} \right]
\int_{-\infty}^{\infty} \exp\left[-\frac{1}{2}
\left(t-\frac{2\sqrt{n}}{\beta}
\right)^2\right] \frac{dt}{\sqrt{2\pi}}\label{pdfInt}\\
&=\frac{\beta}{2\sqrt{n}} \exp\left[\frac{2n}{\beta^2} \right],\nonumber
\end{align}
which gives
\begin{align*}
2\pi \int_0^{\infty}rG \left[\frac{s-\beta\log(Kr)-n/\beta}{\sqrt{n}}
\right ]dr &=  \frac{\pi }{K^2} e^{\frac{2s}{\beta} }.
\end{align*}
We proceed similarly with the other integral in Lemma~\ref{Lemma1}
\begin{align*}
&\int_{B_0(b_n)\setminus B_0(a_n)}\nu_n (s,|x|)dx \\
&= 2\pi
\int_{a_n}^{b_n}rG \left[\frac{s-\beta\log(Kr)-n/\beta}{\sqrt{n}} \right]dr.
\end{align*}
The change of variables $t=\frac{s-\beta\log(Kr)-n/\beta}{\sqrt{n}}$ gives 
\begin{align*}
&2\pi
\int_{a_n}^{b_n}rG \left[\frac{s-\beta\log(Kr)-n/\beta}{\sqrt{n}} \right ]dr\\
&=2\pi \int_{v_n}^{u_n} \frac{1}{K^2}\exp\left[
\frac{2}{\beta}\left(s-t\sqrt{n} -n/\beta \right) \right] 
G (t)\frac{\sqrt{n}}{\beta}dt \\
&=2\pi\frac{\sqrt{n}}{\beta} \frac{\exp\left[ \frac{2}{\beta}\left(s
-n/\beta \right) \right]}{K^2}
 \int_{v_n}^{u_n} \exp\left[ \frac{-2t\sqrt{n}}{\beta} \right]  G (t)dt 
\end{align*}
where 
\begin{align*}
&u_n=\frac{s-\beta\log(K a_n )-n/\beta}{\sqrt{n}}, \\
&v_n=\frac{s-\beta\log(Kb_n)-n/\beta}{\sqrt{n}}.
\end{align*}
Moreover
\begin{align*}
 &\int_{v_n}^{u_n} \exp\left[ \frac{-2t\sqrt{n}}{\beta} \right]  G (t)dt\\
 =&
 \int_{-u_n}^{-v_n} \exp\left[ \frac{2t\sqrt{n}}{\beta} \right] 
\bar{G }(t)dt\\ 
=&\left.\frac{\beta}{2\sqrt{n}} \exp\left[ \frac{2t\sqrt{n}}{\beta} \right] 
\bar{G }(t)\right| _{-u_n}^{-v_n} 
\!\!\!\!\!+\frac{\beta}{2\sqrt{n}}\int_{-u_n}^{-v_n}\exp\left[ \frac{2t\sqrt{n}}{\beta}
\right]  g (t)dt\\
=&\left.\frac{\beta}{2\sqrt{n}} \exp\left[ \frac{2t\sqrt{n}}{\beta} \right] 
\bar{G }(t)\right| _{-u_n}^{-v_n} \\
&+\frac{\beta}{2\sqrt{n}} \exp\left[\frac{2n}{\beta^2} \right]
 \int_{-u_n}^{-v_n} \exp\left[-\frac{1}{2} \left(t-\frac{2\sqrt{n}}{\beta}
\right)^2\right] \frac{dt}{\sqrt{2\pi}}.
\end{align*}
We now derive the conditions of $a_n$ and $b_n$ which ensure that the above integral converges and the integrated term disappears. The latter can be achieved, given inequality (\ref{errIneq}), in the limit as $-u_n$ and $-v_n$ both approach infinity, or equivalently 
$\frac{\beta\log(K b_n)}{\sqrt{n}} +\frac{\sqrt{n}}{\beta}\rightarrow
\infty \quad$ and 
$\frac{\beta\log(K a_n)}{\sqrt{n}} +\frac{\sqrt{n}}{\beta}\rightarrow
\infty$
as $n\rightarrow \infty$, which agrees with conditions (\ref{a_n}) and (\ref{b_n}). Further conditions are revealed after the change of variable $w=t-2 \sqrt{n}/\beta$, yielding the integral
$$
 \int_{-u_n}^{-v_n} e^{\left.-\frac{1}{2} \left(t-\frac{2\sqrt{n}}{\beta}
\right)^2\right.} \frac{dt}{\sqrt{2\pi}}= \int_{-u_n-\frac{2\sqrt{n}}{\beta}}^{-v_n-\frac{2\sqrt{n}}{\beta}} 
e^{-\frac{w^2}{2}} \frac{dw}{\sqrt{2\pi}}
$$
whose limits of integration imply, in light of the earlier integral (\ref{pdfInt}), that as $n\rightarrow \infty$ the following
$\frac{\beta\log(K b_n)}{n} >1/\beta$ and 
$\frac{\beta\log(K a_n)}{n} <1/\beta$
are required, of which both conditions (\ref{a_n}) and (\ref{b_n}) satisfy.

\end{IEEEproof}

\begin{lemma}\label{Lemma2}
Assume (\ref{lambdaDef}), (\ref{a_n}) and (\ref{b_n}), then
\begin{align}\label{Lemma2eq1}
&\lim_{n\rightarrow \infty} \sum_{X_i\in \phi\cap (B_0(b_n)\setminus B_0(a_n))}
\nu_n (s,|X_i|) \\
&=\lim_{n\rightarrow \infty} \sum_{X_i\in \phi} \nu_n
(s,|X_i|)=
\Lambda_{\log}((-\infty,s]).\label{Lemma2eq2}
\end{align}
\end{lemma}

\begin{IEEEproof}
For $k\geq 0$ and a fixed $\epsilon>0$, let $r_k=e^{\epsilon k}$ and
$A_k=B_0(r_{k+1})\setminus B_0(r_k)$, and write the summation
in~(\ref{Lemma2eq2}) as  
\begin{align}\nonumber
&\sum_{X_i\in \phi} \nu_n (s,|X_i|)\\
&=\sum_{X_i\in \phi\cap B_0(r_{k_0})} \nu_n
(s,|X_i|) + 
\sum_{k=k_0}^{\infty}\sum_{X_i\in \phi\cap A_k} \nu_n (s,|X_i|) ,
\label{e.2summands}
\end{align}
for some $k_0 \geq 0$, whose value will be fixed later on. In the
limit of $n\rightarrow\infty$, the first summation
in~(\ref{e.2summands}) disappears
\begin{align*}
&\sum_{X_i\in \phi\cap B_0(r_{k_0})} \nu_n (s,|X_i|)  \\
 &=  \sum_{X_i\in \phi\cap
B_0(r_{k_0})}\Prob \left[Z \leq 
\frac{s-\beta\log(K|X_i|)-n/\beta}{\sqrt{n}} \right ]\\
&=  \sum_{X_i\in \phi\cap B_0(r_{k_0})}
G \left[\frac{s-\beta\log(K|X_i|)-n/\beta}{\sqrt{n}}
\right ]\\
& \leq \phi(B_0(r_{k_0}))G \left[\frac{s-\beta\log(K|X_*|)-n/\beta}{\sqrt{n}} \right ]
\\
&\rightarrow 0 \qquad (n \rightarrow \infty).
\end{align*}
where $X_*$ gives the maximum of $G
\left[\frac{s-\beta\log(K|X|)-n/\beta}{\sqrt{n}} \right ]$ over
$X\in\phi\cap B_0(r_{k_0})$  which exists since $\phi$ is (by
our assumption) a locally finite point measure. 
For the second summation in~(\ref{e.2summands}) we write $\nu_n (s,|X_i|)=\nu_n
(s,|x|\frac{|X_i|}{|x|})$, hence
\begin{align*}
 &\nu_n (s,|X_i|)\\
&=\frac{1}{|A_k|}\int_{A_k}\nu_n (s,|x|\frac{|X_i|}{|x|}) dx.
\end{align*}
Then the bounds
$$
e^{-\epsilon}=\frac{r_{k}}{r_{k+1}}\leq\frac{|X_i|}{|x|} \leq
\frac{r_{k+1}}{r_{k}}=e^{\epsilon},
$$
and form of $\nu_n$, which implies $\nu_n (s,|x|e^{\epsilon})=\nu_n
(s-\beta\epsilon,|x|)$, lead to the lower bound
\begin{equation}
\sum_{k=k_0}^{\infty}\sum_{X_i\in \phi\cap A_k} \nu_n (s,|X_i|) 
\geq
\sum_{k=k_0}^{\infty} \frac{\phi( A_k)}{|A_k|}\int_{A_k}\nu_n
(s-\beta\epsilon ,|x|) dx,
\label{lowerBound1}   
\end{equation}
and the upper bound
\begin{equation}
\sum_{k=k_0}^{\infty}\sum_{X_i\in \phi\cap A_k} \nu_n (s,|X_i|)  \\
\leq  \sum_{k=k_0}^{\infty} \frac{\phi( A_k)}{|A_k|}\int_{A_k}\nu_n
(s+\beta\epsilon,|x|) dx. \label{upperBound1}
\end{equation}
Moreover, we write
\begin{align*} 
\frac{\phi( A_k)}{|A_k|} &= \frac{\phi( B_0(r_{k+1}))-\phi(
B_0(r_{k}))}{|B_0(r_{k+1})|-|B_0(r_{k})|}\\
&= \frac{\frac{\phi( B_0(r_{k+1}))}{|B_0(r_{k+1})|}-\frac{\phi(
B_0(r_{k}))}{|B_0(r_{k})|}\frac{|B_0(r_{k})|}{|B_0(r_{k+1})|}}{1-\frac{|B_0(r_{k
})|}{|B_0(r_{k+1})|}}\\
&= \frac{\frac{\phi( B_0(r_{k+1}))}{|B_0(r_{k+1})|}-\frac{\phi(
B_0(r_{k}))}{|B_0(r_{k})|}e^{-2\epsilon}}{1-e^{-2\epsilon}},
\end{align*}
and requirement (\ref{lambdaDef}) yields
$
\lim_{k\rightarrow \infty}\frac{\phi( A_k)}{|A_k|} = \lambda
$.
Hence, for any fixed $\delta>0$, there exists a $k_0(\delta)$ such that for
all $k\geq k_0$, the bounds
$$
(1-\delta)\lambda \leq \frac{\phi( A_k)}{|A_k|} \leq (1+\delta)\lambda,
$$
hold. Lower bound (\ref{lowerBound1}) becomes
\begin{align*} 
\sum_{k=k_0}^{\infty}\sum_{X_i\in \phi\cap A_k} \!\!\!\nu_n (s,|X_i|)
&\geq
\sum_{k=k_0}^{\infty} (1-\delta)\lambda \int_{A_k}\!\!\!\nu_n (s-\beta\epsilon,|x|)
dx,\\  
&=(1-\delta)\lambda \int_{|x|\geq r_{k_0}}\!\!\!\!\!\nu_n (s-\beta\epsilon,|x|) dx
\end{align*}
Finally, Lemma \ref{Lemma1} allows us to set $a_n=r_{k_0}$ and $b_n=\infty$, hence 
$$
 \lim_{n\rightarrow \infty}\sum_{k=k_0}^{\infty}\sum_{X_i\in \phi\cap A_k}
\nu_n (s,|X_i|)\geq (1-\delta) \frac{\pi\lambda
 }{K^2} e^{\frac{2s-\beta\epsilon}{\beta} },
$$
and similarly the upper bound (\ref{upperBound1}) becomes
$$ 
\lim_{n\rightarrow \infty}\sum_{k=k_0}^{\infty}\sum_{X_i\in \phi\cap A_k}
\nu_n (s,|X_i|)\leq (1-\delta) \frac{\pi\lambda
 }{K^2} e^{\frac{2s+\beta\epsilon}{\beta} },
$$
and indeed $\epsilon\rightarrow 0$ and $\delta \rightarrow 0$
completes the proof of (\ref{Lemma2eq2}).
The other result, (\ref{Lemma2eq1}), can be proved by a
straightforward modification
of the above arguments.  
\end{IEEEproof}

\begin{IEEEproof}[Proof of Theorem \ref{mainResult}]
We use a classical convergence result~\cite[Theorem 11.22V]{daleyPPII2008} 
which in our  setting requires verification of
the following two conditions
 (cf~\cite[(11.4.2) and (11.4.3)]{daleyPPII2008})
\begin{equation}\label{cond1}
\sup_i \nu^i_n (A)  \rightarrow 0\quad (n \rightarrow \infty).
\end{equation}
and
\begin{equation}\label{cond2}
\sum_i \nu^i_n (A) \rightarrow \Lambda_{\log}(A) \quad
(n \rightarrow \infty),
\end{equation}
for all bounded Borel sets $A\subset\R$, 
where $\nu_n^i(\cdot)$  is the (probability) measure on $\R$ defined
by setting $\nu_n^i((-\infty,s]):=\nu_n(s,|X_i|)$.
The first condition, (\ref{cond1}), clearly holds by~(\ref{v_nFinal}) for any locally finite
$\phi$. 
The second condition, (\ref{cond2})
follows from Lemma \ref{Lemma1} and \ref{Lemma2}, which establish the
required convergence for $A=(-\infty,s]$ and any~$s\in\R$. This is
enough to conclude the convergence for all bounded Borel sets.
\end{IEEEproof}

\subsection{Representation of the function~$\varphi_{\beta}(z)$ given by~(\ref{e.LaplaceF0})}
\label{App2}
\begin{lemma}
The function $\varphi_{\beta}\left(  z\right)  $\ given by~(\ref{e.LaplaceF0})
can be written as 
\begin{align}
\varphi_{\beta}\left(  z\right)   &  =-\frac{2}{\beta}z^{\frac{2}{\beta}%
}\gamma\left(  -\frac{2}{\beta},z\right) \label{e.Phi1}\\
&  =\Gamma\left(  1-\frac{2}{\beta}\right)  \gamma^{\ast}\left(  -\frac
{2}{\beta},z\right)  \label{e.Phi2}%
\end{align}
where
\begin{align}
\gamma^{\ast}\left(  \alpha,z\right)   &  =\frac{z^{-\alpha}}{\Gamma\left(
\alpha\right)  }\gamma\left(  \alpha,z\right) \label{e.ModifiedGamma1}\\
&  =e^{-z}\sum_{k=0}^{\infty}\frac{z^{k}}{\Gamma\left(  \alpha+k+1\right)  }
\label{e.ModifiedGamma2}%
\end{align}
is an analytic function on $\mathbb{C}\times\mathbb{C}$. Moreover
\begin{equation}
\varphi_{\beta}\left(  z\right)  >0,\quad\forall\beta>2,\forall z\in
\mathbb{R}^{+} \label{e.PhiPositive}%
\end{equation}
\end{lemma}

\begin{IEEEproof}
Using~(\ref{e.LaplaceF0}) and recurrence formula~\cite[6.5.22]%
{Abramowitz1970}\ we obtain~(\ref{e.Phi1}). Introducing the modified incomplete
gamma function $\gamma^{\ast}$\ defined by~(\ref{e.ModifiedGamma1}) (which is
holomorphic on $\mathbb{C}\times\mathbb{C}$; see~\cite[6.5.4]{Abramowitz1970})
we obtain~(\ref{e.Phi2}). The expansion~(\ref{e.ModifiedGamma2}) is given
in~\cite[6.5.29]{Abramowitz1970}. This expansion shows that if $\alpha>-1$
then for any $z\in\mathbb{R}^{+}$, $\gamma^{\ast}\left(  \alpha,z\right)  >0$,
and therefore by~(\ref{e.Phi2}) $\varphi_{\beta}\left(  z\right)  >0$ for any
$\beta>2,z\in\mathbb{R}^{+}$.
\end{IEEEproof}

\pdfbookmark[0]{References}{References}


\end{document}